\documentstyle[11pt]{article}

\author{A. D. King\thanks{Supported by a research grant from the SERC.}\\
Dept. of Pure Mathematics\\
University of Liverpool\\
P. O. Box 147\\
Liverpool L69 3BX\\
UNITED KINGDOM\\
e-mail: aking@liverpool.ac.uk
\and
Charles H. Walter\thanks{Supported in part by NSA research grant
MDA904-92-H-3009.}\\
URA 168\\
Math\'ematiques\\
Universit\'e de Nice\\
F-06108 Nice Cedex 02 \\
FRANCE\\
e-mail: walter@math.unice.fr}
\title{On Chow Rings of Fine Moduli Spaces of Modules
}
\date{
}
\newtheorem{theorem}{Theorem}

\newcommand{\Rsym}{{\cal R}^{\rm sym}(Q,\alpha)}
\newcommand{\Msym}{M_{S}(\alpha,\theta)}
\def\qed{\hfill\mbox{$QED$}}
\def\qed{\hfill\mbox{$\Box$}}
\input{amssym.def}
\newsymbol\ncong 231D
\def\lrabox
{\hbox{$\longrightarrow$}}
\def\triplearrow
{\vcenter{\lrabox\nointerlineskip\lrabox\nointerlineskip\lrabox}}
\def\doublearrow
{\vcenter{\lrabox\nointerlineskip\lrabox}}
\catcode`\@=11
\def\@begintheorem#1#2{\sl \trivlist
   \item[\hskip \labelsep{\bf #1\ #2\thmcounterend}]}
\def\@opargbegintheorem#1#2#3{\sl \trivlist
      \item[\hskip \labelsep{\bf #1\ #2\ (#3)\thmcounterend}]}
\def\thmcounterend{.}
\def\section{\@startsection{section}{1}{\z@}{-3.25ex plus
 -1ex minus -.2ex}{1.5ex plus .2ex}{\large\bf}}
\def\subsection{\@startsection
 {subsection}{2}{\z@}{3.25ex plus 1ex minus .2ex}{-0.5em}{\normalsize\sl}}
\def\subsubsection{\@startsection
 {subsubsection}{3}{\z@}{3.25ex plus 1ex minus .2ex}{-0.5em}{\normalsize\sl}}
\def\paragraph{\@startsection
 {paragraph}{3}{\z@}{2ex plus 0.6ex minus .2ex}{-0.5em}{\normalsize\sl}}
\def\subparagraph{\@startsection
 {subparagraph}{3}{\parindent}{2ex plus 0.6ex minus
     .2ex}{-1pt}{\normalsize\sl}}
\catcode`\@=12

\begin{document}

\maketitle

\begin{abstract}
Let $M$ be a complete nonsingular fine moduli space of modules over an
algebra $S$.  A set of conditions is given for the Chow ring of $M$ to
be generated by the Chern classes of certain universal bundles
occurring in a projective resolution of the universal $S$-module on
$M$.  This result is then applied to the varieties $G_T$ parametrizing
homogeneous ideals of $k[x,y]$ of Hilbert function $T$, to moduli
spaces of representations of quivers, and finally to moduli spaces of
sheaves on ${\Bbb P}^2$, reinterpreting a result of Ellingsrud and
Str\o mme.
\end{abstract}

In a recent paper \cite{ES} Ellingsrud and Str\o mme identified a set
of generators of the Chow ring of the moduli space of stable sheaves
of given rank and Chern classes on ${\Bbb P}^2$ (in the case where the
moduli space is smooth and projective). In this paper we formulate a
part of their argument as a general theorem about fine moduli spaces
of modules over an associative algebra.  This provides a more widely
applicable method for showing that the Chow ring of a fine moduli
space is generated by the Chern classes of appropriate universal
sheaves.  In particular we apply the method to verify a conjecture of
Iarrobino and Yam\'eogo concerning the Chow rings of the varieties
$G_T$ parametrizing homogeneous ideals in $k[x,y]$ with a given
Hilbert function. We also verify a conjecture of the first author
concerning Chow rings of moduli spaces of representations of quivers.

Let $S$ be an associative algebra over an algebraically closed field
$k$.  By convention we will consider only left $S$-modules in this
paper.  A flat family of $S$-modules over a $k$-scheme $X$ is a sheaf
$\cal F$ of $S\otimes {\cal O}_X$-modules on $X$, quasi-coherent and
flat over ${\cal O}_X$.  At a (closed) point $x\in X$ the fiber ${\cal
F}(x)$ is an $S$-module.  If $\cal C$ is a class of $S$-modules, then
a fine moduli space for $\cal C$ is a scheme $M$ equipped with a flat
family $\cal U$ all of whose fibers are in $\cal C$ and with the usual
universal property.  Our general theorem is the following:

\begin{theorem}
\label{main}Let ${\cal C}$ be a class of $S$-modules, and $M$ a fine
moduli space for ${\cal C}$ which is a complete nonsingular
variety.  Suppose further that

(i)\quad If $E\in {\cal C}$, then ${\rm Hom}_S(E,E)\cong k$, ${\rm Ext}
_S^1(E,E)\cong T_{[E]}M$, and ${\rm Ext}_S^p(E,E)=0$ for $p\geq 2$;

(ii)\quad If $E\ncong F$ are in ${\cal C}$, then ${\rm Hom}_S(E,F)=0$
and ${\rm Ext}_S^p(E,F)=0$ for $p\geq 2.$

(iii)\quad If ${\cal U}$ is the universal $S\otimes {\cal O}_M$-module
on $M$, then ${\cal U}$ has a universal projective resolution of
finite length $$ 0\rightarrow \bigoplus_jP_{rj}\otimes {\cal
E}_{rj}\rightarrow \cdots
\rightarrow \bigoplus_jP_{1j}\otimes {\cal E}_{1j}\rightarrow
\bigoplus_jP_{0j}\otimes {\cal E}_{0j}\rightarrow {\cal U}\rightarrow 0
$$
with the $P_{ij}$ projective $S$-modules such that $\dim _k{\rm Hom}
_S(P_{ij},P_{i^{\prime }j^{\prime }})$ is always finite, and the ${\cal E}
_{ij}$ are all locally free ${\cal O}_M$-modules of finite rank.

Then

(a)\quad The Chern classes of the ${\cal E}_{ij}$ generate the Chow ring $
A^{*}(M)$ as a ${\Bbb Z}$-algebra.

(b)\quad Numerical and rational equivalence coincide on $M$. In particular, $
A^{*}(M)$ is a free ${\Bbb Z}$-module.

(c)\quad If $k={\Bbb C}$, then the cycle map $A^{*}(M)\rightarrow H^{*}(M,{
\Bbb Z})$ is an isomorphism. In particular, there is no odd-dimensional
cohomology.
\end{theorem}

If $S$ is a graded algebra, then one may formulate a graded version of
the theorem by replacing ``module'' by ``graded module'' throughout
and using the degree-zero parts of ${\rm Hom}$ and ${\rm Ext}$.  More
generally, one should be able to formulate a version of Theorem
\ref{main} for moduli spaces of objects in an abelian category of
$k$-vector spaces in any situation where one has a suitable notion of
a family of objects of the category.

Our first application of Theorem \ref{main} is to the Iarrobino
varieties $G_T$ which parametrize homogeneous ideals $I\subset k[x,y]$
of Hilbert function $T$.  Here $T=(t_0,t_1,\ldots )$ is a sequence of
nonnegative integers such that $t_n=0$ for $n\gg 0$, and the points of
$G_T$ correspond to those $I$ such that $\dim_k(k[x,y]/I)_n = t_n$ for
all $n$. These smooth projective varieties were originally constructed
in \cite{I} in order to study the Hilbert-Samuel function
stratification of the punctual Hilbert scheme of a surface. They have
since been studied in several papers including \cite{Goe} \cite{IY}
and \cite{Y}. The fact that these $G_T$ are fine moduli spaces was
addressed formally in \cite{G} Kap.\ 2, Lemma 4.

Having fixed the Hilbert function $T$, the degree-$n$ graded pieces of
the quotient rings form a family of quotients of dimension $ t_n$ of
the space of binary forms of degree $n$. This induces a natural
morphism from $G_T$ to the Grassmannian of quotient spaces ${\rm Gr}
(t_n,n+1) $. Let ${\cal A}_n$ denote the pullback to $G_T$ of the
universal quotient bundle on ${\rm Gr}(t_n,n+1)$. Our result is

\begin{theorem}
\label{GT}Let $G_T$ be the Iarrobino variety parametrizing homogeneous
ideals $I\subset k[x,y]$ of Hilbert function $T$, and let the ${\cal A}_n$
be the universal bundles defined above. Then

(a)\quad The Chern classes of the ${\cal A}_n$ generate the Chow ring $
A^{*}(G_T)$ as a ${\Bbb Z}$-algebra.

(b)\quad Numerical and rational equivalence coincide on $G_T$. In
particular, $A^{*}(G_T)$ is a free ${\Bbb Z}$-module.

(c)\quad If $k={\Bbb C}$, then the cycle map $A^{*}(G_T)\rightarrow
H^{*}(G_T,{\Bbb Z})$ is an isomorphism. In particular, there is no
odd-dimensional cohomology.
\end{theorem}

Parts (b) and (c) were already known because $G_T$ has a cell decomposition
corresponding to the initial ideals with Hilbert function $T$ (cf.\
\cite{Goe} or \cite{IY}). Nevertheless, our methods give a new proof.

Our second application is to fine moduli spaces of representations of
a quiver without oriented cycles.  These moduli spaces were
constructed in \cite{K}.

To fix notation, we recall that a quiver $Q$ is a directed graph,
specified by a finite set of vertices $Q_0$ and a finite set of arrows
$Q_1$ between the vertices together with two maps
$h,t{:}~Q_1\rightarrow Q_0$ specifying the head and tail of each
arrow.  A representation of $Q$ consists of vector spaces $W_i$ for
each $i\in Q_0$ and $k$-linear maps $\phi_a:W_{ta}\to W_{ha}$ for each
$a\in Q_1$.  A subrepresentation is a collection of subspaces
$W_i'\subset W_i$ such that $\phi_a(W_{ta}')\subset W_{ha}'$ for all $a$.
The dimension vector $\alpha \in {\Bbb N}^{Q_0}$ of a representation
$(W_i,\phi_a)$ is given by $\alpha_i = \dim_k(W_i)$.

To obtain a moduli space of representations of $Q$ of dimension vector
$\alpha$ one needs to introduce a notion of stability.  Having chosen
$\theta = (\theta_i)\in {\Bbb R}^{Q_0}$ such that $\sum_i
\theta_i\alpha_i = 0$, we say that a representation $(W_i,\phi_a)$ is
$\theta$-stable if all (proper) subrepresentations $(W_i')$ satisfy
$\sum_i \theta_i \dim(W_i) > 0$.  When $\alpha $ is an indivisible
dimension vector and $\theta $ is generic, there is a smooth fine
moduli space $M_Q(\alpha,\theta)$ of $\theta $-stable representations
of $Q$ of dimension vector $\alpha $ (\cite{K} Proposition 5.3).  If
the quiver $Q$ has no oriented cycles, then this fine moduli space is
projective (\cite{K} Proposition 4.3).  Note that this moduli
space may actually be empty.  The conditions on $\alpha$ and $\theta$
which make it non-empty are more subtle (cf.\ \cite{K} Remark 4.5).

The universal representation over $M_Q(\alpha ,\theta )$ consists of
vector bundles ${\cal U}_i$ of rank $\alpha_i$ together with the
universal morphisms.  We use Theorem \ref{main} to prove the
following, confirming the conjecture made in Remark 5.4 of \cite{K}.

\begin{theorem}\label{quiver}
Let $Q$ be a quiver without oriented cycles, and
$M=M_Q(\alpha,\theta)$ be a smooth projective fine moduli space of
$\theta$-stable representations of $Q$ of dimension vector $\alpha$.
Let ${\cal U}_i$ be the universal bundles on $M$ described above.  Then

(a)\quad The Chern classes of the ${\cal U}_i$ generate the Chow ring
$A^{*}(M)$ as a ${\Bbb Z}$-algebra.

(b)\quad Numerical and rational equivalence coincide on $M$.
In particular, $A^{*}(M)$ is a free ${\Bbb Z}$-module.

(c)\quad If $k={\Bbb C}$, then the cycle map $A^{*}(M)\to H^{*}(M,{\Bbb Z})$
is an isomorphism. In particular, there is no odd-dimensional cohomology.
\end{theorem}

The outline of the paper is as follows. In the first section we prove
Theorem \ref{main} by adapting the method of Ellingsrud and Str\o mme. In
the second and third sections we apply Theorem \ref{main} to prove Theorems
\ref{GT} and \ref{quiver}. In the fourth section we explain how
Theorem \ref{main} may be used to prove Ellingsrud and Str\o mme's
original result for sheaves on ${\Bbb P}^2$.

\section{Proof of the Main Theorem}

In this section we prove Theorem \ref{main} by adapting a method of
Ellingsrud and Str\o mme.

Let $\delta $ be the class of the diagonal in $A^{*}(M\times M)$, and let $
p_1$ and $p_2$ denote the projections from $M\times M$ onto its two factors.
We will adapt the methods of \cite{ES} \S 2 to show that $\delta $ can be
written as a polynomial in the Chern classes of the $p_1^{*}({\cal E}_{ij})$
and the $p_2^{*}({\cal E}_{ij})$. The theorem will then follow from the
following result, which Ellingsrud and Str\o mme describe as ``a well-known
observation on varieties with decomposable diagonal class'':

\begin{theorem}
{\rm (\cite{ES} Theorem 2.1)} Let $X$ be a nonsingular complete variety.
Assume that the rational equivalence class $\delta $ of the diagonal $\Delta
\subseteq X\times X$ decomposes in the form
\begin{equation}
\label{diag}\delta =\sum_{i\in J}p_1^{*}\alpha _i\ p_2^{*}\beta _i
\end{equation}
where $p_1$ and $p_2$ are the projection of $X\times X$ onto its factors,
and $\alpha _i$, $\beta _i\in A^{*}(X)$. Then

(a)\quad The $\alpha _i$ generate $A^{*}(X)$ as a ${\Bbb Z}$-module.

(b)\quad Numerical and rational equivalence coincide on $X$. In particular, $
A^{*}(X)$ is a free ${\Bbb Z}$-module.

(c)\quad If $k={\Bbb C}$, then the cycle map $A^{*}(X)\rightarrow H^{*}(X,{
\Bbb Z})$ is an isomorphism. In particular, there is no odd-dimensional
cohomology.

(d)\quad Suppose the set $\{\alpha _i\}$ in (\ref{diag}) is minimal. Then $
\{\alpha _i\}$ and $\{\beta _i\}$ are dual bases with respect to the
intersection form on $A^{*}(X)$.
\end{theorem}

So we show how to write $\delta $ as a polynomial in the Chern classes
of the $p_1^{*}({\cal E}_{ij})$ and the $p_2^{*}({\cal E}_{ij})$ using
a method similar to \cite{ES} \S 2. First we write ${\cal P}_{\bullet
}$ for the projective resolution of the universal family of modules
${\cal U}$ in unaugmented form
$$
0\rightarrow \bigoplus_jP_{rj}\otimes {\cal E}
_{rj}\rightarrow \cdots \rightarrow \bigoplus_jP_{0j}\otimes {\cal E}
_{0j}\rightarrow 0.
$$
Then let ${\cal L}^{\bullet}={\cal H}om_{S\otimes {\cal O}_{M\times
M}}^{\bullet}(p_1^{*}{\cal P}_{\bullet},p_2^{*}{\cal P}_{\bullet})$.
Since
\begin{equation}
\label{Lp}{\cal L}^p=\bigoplus_{i^{\prime }-i=p}\left(
\bigoplus_{j,j^{\prime }}{\rm Hom}_S(P_{ij},P_{i^{\prime }j^{\prime
}})\otimes p_1^{*}{\cal E}_{ij}^{\vee }\otimes p_2^{*}{\cal
E}_{i^{\prime }j^{\prime }}\right) ,\
\end{equation}
${\cal L}^{\bullet}$ is a finite complex of locally free modules of
finite rank with the universal property that for any morphism of
$k$-schemes of the form $\phi {:}\ X\rightarrow M\times M$, we have
$H^p(\phi ^{*}{\cal L} ^{\bullet})={\cal E}xt_{S\otimes {\cal
O}_X}^p(\phi ^{*}p_1^{*}{\cal U},\phi ^{*}p_2^{*}{\cal U})$ for all
$p$. In particular ${\cal L}^{\bullet}$ is exact except in degrees 0
and 1. Indeed, if $d^p{:}\ {\cal L}^p\rightarrow {\cal L}^{p+1}$ is
the differential of ${\cal L}^{\bullet}$, then ${\cal L} ^{\bullet}$
is quasi-isomorphic to the short complex
$$
0\rightarrow {\rm cok}(d^{-1})\stackrel{\phi }{\longrightarrow }{\rm
\ker } (d^1)\rightarrow 0
$$
where $\phi $ is a map between locally free sheaves whose degeneracy
locus is exactly the diagonal. Our complex ${\cal L}^{\bullet}$ now
has all the essential properties of the complex ${\cal C}^{\bullet}$
of \cite{ES} Lemma 2.4. Hence by the same argument we have
$$
\delta =c_{\dim M}([{\rm \ker }(d^1)]-[{\rm cok}(d^{-1})])=c_{\dim
M}(\sum (-1)^{p+1}[{\cal L}^p]).  $$ The formula (\ref{Lp}) and
standard formulas for Chern classes now permit us to write $\delta $
as a polynomial in the Chern classes of the $p_1^{*}{\cal E}_{ij}$ and
the $p_2^{*}{\cal E}_{i^{\prime }j^{\prime }}$. This completes the
proof of Theorem \ref{main}. {\qed}

\section{Iarrobino Varieties\label{iarro}}

In this section we apply the main theorem to the study of the
Iarrobino varieties $G_T$ which parametrize homogeneous ideals in
$k[x,y]$ with Hilbert function $T$. We prove Theorem \ref{GT}.

\paragraph{Proof of Theorem \ref{GT}.}

Let $S=k[x,y]$ and ${\cal S}=S\otimes _k{\cal O}_{G_T}$. The fact that
$G_T$ is a fine moduli space means that there is a universal sheaf of
homogeneous ${\cal S}$-ideals ${\cal I}$ and a universal quotient
sheaf of graded ${\cal S}$-modules ${\cal A}={\cal S}/{\cal
I}$. According to our construction, we have ${\cal A}=\bigoplus_n{\cal
A}_n$.

We now wish to apply Theorem \ref{main} using ${\cal A}$ as the
universal module. To do so we need to verify the various hypotheses.

First the cohomological ones. Because we are working with graded
modules, we must examine the degree 0 graded pieces of the internal
Hom and Ext. If $ A=S/I$ and $B=S/J$ are two graded quotients of $S$
with the same Hilbert function, then there are no nonzero morphisms of
degree 0 between $A$ and $B$ unless $A=B$, in which case ${\rm
Hom}_S(A,A)_0\cong k$, the homotheties. Some standard exact sequences
show that ${\rm Ext}_S^1(A,A)_0\cong {\rm Hom} _S(I,A)_0$ which is the
tangent space to $G_T$ at $[I]$ by the graded analog of Grothendieck's
formula for the tangent space of the Hilbert scheme (cf.\ \cite{PS} \S
4).  Since $S$ is of global dimension $2$, the functor ${\rm Ext}^2_S$
is right exact, so the surjection $S \to B$ induces a surjection ${\rm
Ext}^2_S(A,S)_0 \to {\rm Ext}^2_S(A,B)_0$.  But by local duality
(\cite{S} n$^{\circ}$ 72, Th\'eor\`eme 1) we have ${\rm
Ext}^2_S(A,S)_0^* \cong H^0_{\frak m}(A)_{-2} = A_{-2}$, which
vanishes.  So ${\rm Ext}^2_S(A,B)_0 = 0$. Finally ${\rm
Ext}_S^p(A,B)=0$ for all $p\geq 3$ because $S$ has global dimension
$2$. Thus the cohomological hypotheses of Theorem \ref {main} are
fulfilled.

Now we exhibit a universal projective resolution of ${\cal A}$. Since
${\cal A}=\bigoplus {\cal A}_n$ is a $k[x,y]\otimes {\cal
O}_{G_T}$-module, multiplication by $x$ and $y$ define morphisms $\xi
$, $\eta {:}\ {\cal A} _n\rightarrow {\cal A}_{n+1}$. Then the
universal projective resolutions is
$$
0\rightarrow \bigoplus_nS(-n-2)\otimes _k{\cal A}_n\stackrel{\alpha }{
\longrightarrow }\bigoplus_nS(-n-1)^2\otimes _k{\cal
A}_n\stackrel{\beta }{\longrightarrow }\bigoplus_nS(-n)\otimes _k{\cal
A}_n\rightarrow {\cal A}\rightarrow 0
$$
where the morphisms are the standard ones
$$
\alpha =\left[
\begin{array}{c}
-y\otimes 1+1\otimes \eta \\
x\otimes 1-1\otimes \xi
\end{array}
\right] ,\qquad \beta =\left[
\begin{array}{cc}
x\otimes 1-1\otimes \xi & y\otimes 1-1\otimes \eta
\end{array}
\right] .
$$
Note that the sums are finite because each ${\cal A}_n$ is of rank $t_n$
which vanishes for $n\gg 0$. Finally the ${\rm Hom}_S(S(-i),S(-j))_0$ are
all finite-dimensional. So the remaining hypotheses of Theorem \ref{main}
are fulfilled.

The theorem now follows directly from Theorem \ref{main}. {\qed}

\section{Representations of Quivers\label{reps}}

In this section we prove Theorem \ref{quiver} by applying the main
theorem to fine moduli spaces of representations of quivers without
oriented cycles.  To do this, we start with the well-established
observation (cf.\ for example \cite{Ben}) that representations of a
quiver $Q$ are the same as modules for the path algebra
$kQ$.  This algebra is generated over $k$ by a set of orthogonal
idempotents $\{e_i\mid i\in Q_0\}$ and a further set of generators
$\{x_a\mid a\in Q_1\}$ such that $x_a=e_{ha}x_{a}e_{ta}$.  A left
$kQ$-module $E$ corresponds to the representation of $Q$ consisting of
the vector spaces $W_i=e_iE$ for each $i\in Q_0$, and the $k$-linear
maps $\phi_a:W_{ta}\to W_{ha}$ giving multiplication by $x_a$ for each
$a\in Q_1$.

The algebra $kQ$ is always hereditary, i.e.\ of global
dimension $\leq1$, and is finite-dimensional if and only if the quiver
$Q$ has no oriented cycles.

\paragraph{Proof of Theorem \ref{quiver}.}

Let $S=kQ$.  Then the indecomposable projective $S$-modules are $Se_i$
for $i\in Q_0$, and $S$ has the following minimal projective
resolution as an $S,S$-bimodule (or $S\otimes S^{op}$-module)
\[
 0\rightarrow \bigoplus_{a\in
Q_1}Se_{ha}\otimes e_{ta}S\;{\frac{\buildrel d}{
\longrightarrow }}\;\bigoplus_{i\in Q_0}Se_i\otimes e_iS\;
{\frac{\buildrel \mu }{\longrightarrow }}\;S\rightarrow 0
\]
where $\mu $ is multiplication and $d(e_{ha}\otimes e_{ta})=x_a\otimes
e_{ta}-e_{ha}\otimes x_a$.  We can use this resolution to calculate
the derived functors of ${\rm Hom}_S$, because ${\rm Hom}_S(E,F)={\rm
Hom}_{S,S}(S,{\rm Hom}_k(E,F))$ for any $E$ and $F$.  We see
immediately see that ${\rm Ext}_S^i(E,F)=0$, for $i\geq 2$.

To check the other cohomological conditions, we first note that, by a
standard Schur's Lemma style argument, the stability condition implies
that for any two $\theta$-stable modules $E$ and $F$
$$
{\rm Hom}_S(E,F)=\cases{ k & if $E\cong F$, \cr 0 & otherwise.}
$$
Now $M$ is constructed as a GIT quotient of the
representation space
$$
{\cal R}(Q,\alpha )=\bigoplus_{a\in Q_1}{\rm Hom}(W_{ta},W_{ha})
$$
by the reductive group $GL(\alpha )=\prod_{i\in Q_0}GL(W_i)$, where
$W_i$ is a fixed vector space of dimension $\alpha _i$.  If
$E=(W_i,\phi_a)$, then the tangent space $T_{[E]}(M)$ is
isomorphic to normal space at $\phi\in {\cal R}(Q,\alpha)$ to the
$GL(\alpha)$-orbit.  This is the cokernel of
$$
d_{\phi}: \quad \bigoplus_{i\in Q_0}{\rm Hom}(W_i,W_i) \longrightarrow
\bigoplus_{a\in Q_1}{\rm Hom}(W_{ta},W_{ha})
$$
where $(d_\phi \gamma )_a=\phi _a\gamma _{ta}-\gamma _{ha}\phi _a$.
But this is exactly the complex which calculates ${\rm
Ext}^{*}_S(E,E)$, using the projective resolution of $S$ above. Thus
$T_{[E]}M\cong {\rm Ext}^1_S(E,E)$.

Finally, we obtain the necessary universal projective resolution
by tensoring ${\cal U}$ with the projective bimodule resolution of $S$,
giving
$$
0\rightarrow \bigoplus_{a\in Q_1}Se_{ha}\otimes {\cal
U}_{ta}\rightarrow \bigoplus_{i\in Q_0}Se_i\otimes {\cal
U}_i\rightarrow {\cal U}\rightarrow 0
$$
This completes the verification of the hypotheses for Theorem
\ref{main} and hence the proof of Theorem \ref{quiver}. {\qed}

\section{Sheaves on ${\Bbb P}^2$}

In \cite{ES} Ellingsrud and Str\o mme proved that the Chow ring of a
moduli space $M$ of stable sheaves on ${\Bbb P}^2$ of fixed rank and
Chern classes is generated by the Chern classes of three bundles on
$M$ in those cases where $M$ is smooth and projective.  We show how
their result can be viewed as an application of our Theorem
\ref{main}.

We use notation derived from a recent paper of Le Potier \cite{L}.
Let $r$, $c_1$, $\chi $, and $m$ be integers such that $-r<c_1\leq
0$, $\chi \leq 0$, $\chi \leq r+2c_1$, and $m\gg 0$. Write $ n=-\chi
+r+c_1$.  We consider representations of the quiver with triple edges
labeled $x_1,y_1,z_1$ and $x_2,y_2,z_2$
\[
\begin{array}{ccccccc}
\alpha : &  & n+c_1 &  & n &  & n-(r+c_1) \\
&  & \bullet & \triplearrow & \bullet & \triplearrow & \bullet \\
\theta : &  & -(r+c_1)m+n &  & (r+2c_1)m-2n+r &  & -c_1m+n
\end{array}
\]
with dimension vector $\alpha$ as marked.  Those representations
satisfying the symmetric relations $x_1y_2=y_1x_2$, $x_1z_2=z_1x_2$,
and $y_1z_2=z_1x_2$ form a closed subvariety $\Rsym$ of the
representation space ${\cal R}(Q,\alpha)$.  Let $I$ be the two-sided
ideal of $kQ$ generated by the symmetric relations above, and let $S=
kQ/I$.  For any $\alpha$ (resp.\ $\theta$) we say that an $S$-module
is of dimension vector $\alpha$ (resp.\ is $\theta$-stable) if it is
so as a $kQ$-module.  Then for any $\theta$ the image of $\Rsym$ in
the moduli space $M_Q(\alpha,\theta)$ is a fine moduli space $\Msym$
of $\theta$-stable $S$-modules of dimension vector $\alpha$.

Le Potier has established a result (\cite{L} Th\'eor\`eme 3.1) which
may be interpreted in this language as saying that for $\alpha$ and
$\theta$ as marked in the diagram above, $\Msym$ is isomorphic to the
moduli space $M_{{\Bbb P} ^2}(r,c_1,\chi )$ of Gieseker-Maruyama
stable sheaves on ${\Bbb P}^2$ of rank $r$, determinant $c_1$ and
Euler characteristic $\chi $. The restrictions to $\Msym$ of the
universal bundles ${\cal U}_i$ on $M_Q(\alpha,\theta)$ may be
identified with the universal bundles $R^1\pi _{*}( {\cal E}(-i))$ on
$M_{{\Bbb P}^2}(r,c_1,\chi )$. (Note that with these conventions, the
vertices of the quiver are labeled $2,1,0$ from left to right.)  Lemma
2.2 of \cite{ES} may be interpreted as saying that the Ext groups for
$\theta$-stable $S$-modules are isomorphic to the Ext groups for the
corresponding stable sheaves on ${\Bbb P}^2$.  Hence the cohomological
hypotheses in our Theorem \ref{main} can be verified for
$\theta$-stable $S$-modules by using properties of stable sheaves on
${\Bbb P}^2$.

Let $\overline{e}_i$ be the images in $S$ of the idempotents $e_i$ of
$kQ$ corresponding to the three vertices $2,1,0\in Q_0$ (cf.\
\S\ref{reps}).  Then the indecomposable projective $S$-modules are
$S\overline{e}_i$.  As in \S\ref{reps} the minimal projective
resolution of $S$ as an $S,S$-bimodule yields a projective resolution
of the universal $S$-module $\cal U$ on $\Msym$ which is now of the
form
\[
0\rightarrow \bigoplus_{\mathop{3\
relations}}S\overline{e}_0\otimes{\cal U}_2 \rightarrow
\bigoplus_{a\in Q_1}S\overline {e}_{ha}\otimes {\cal
U}_{ta}\rightarrow \bigoplus_{i\in Q_0}S\overline {e}_i\otimes {\cal
U}_i\rightarrow {\cal U}\rightarrow 0.
\]
We could therefore apply Theorem \ref{main} to retrieve \cite{ES}
Theorem 1.1.

\paragraph{Remark.}

It is also possible to prove Theorem \ref{GT} by regarding the
Iarrobino varieties $G_T$ as moduli spaces for representations of a
quiver ``with relations'' and thus as moduli spaces for
modules over a finite-dimensional non-commutative algebra.  If
$T=(t_0,t_1,\ldots ,t_q,0,0,\ldots )$, then the algebra will be of the
form $R=kQ/I$ where $Q$ is the quiver
\[
\begin{array}{ccccccccccccc}
\alpha : &  & t_0 &  & t_1 &  & t_2 &  &  &  & t_{q-1} &  & t_q \\
&  & \bullet & \doublearrow & \bullet & \doublearrow & \bullet &
\doublearrow & \cdots & \doublearrow & \bullet & \doublearrow &
\bullet \\
\theta : &  & - &  & + &  & + &  &  &  & + &  & +
\end{array}
\]
with $2q$ edges $x_1,y_1,x_2,y_2,\ldots ,x_q,y_q$, and $I$ is
generated by the relations $x_iy_{i+1} - y_ix_{i+1}$.  There is a
clear correspondence between $R$-modules of dimension vector $\alpha$
(as marked) and graded $k[x,y]$-modules of Hilbert function $T$.  When
$t_0=1$ and the coefficients of $\theta $ have the signs indicated
above, then the $\theta $-stable $R$-modules correspond exactly to
$k[x,y]$-modules which are generated in degree $0$, i.e.\ to modules
isomorphic to $k[x,y]/J$ for some homogeneous ideal $J$ of $k[x,y]$.
Hence $M_R(\alpha,\theta)\cong G_T$.  The universal bundles ${\cal
U}_i$ on $M_R(\alpha,\theta)$ are exactly the universal bundles ${\cal
A}_i$ on $G_T$.


\begin{thebibliography}{Ben}

\bibitem[B]{Ben} {\em D.~J.~Benson}, Representations and Cohomology
I: Basic Representation Theory of Finite Groups and Associative
Algebras, Cambridge Studies in Advanced Mathematics, vol.\ 30,
Cambridge University Press, Cambridge, 1991.

\bibitem[ES]{ES}  {\em G.~Ellingsrud} and {\em S.~A.~Str\o mme},
Towards the Chow ring of the Hilbert scheme of ${\Bbb P}^2$, J.\ reine
angew.\ Math.\ {\bf 441} (1993), 33--44.

\bibitem[G]{Goe}  {\em L.~G\"ottsche}, Betti numbers for the Hilbert
function strata of the punctual Hilbert scheme in two variables,
Manuscripta Math.\ {\bf 66} (1990), 253--259.

\bibitem[Go]{G}  {\em G.~Gotzmann}, Topologische Eigenschaften von
Hilbertfunktion-Strata, Habi\-li\-ta\-tions\-schrift, M\"unster, 1993.

\bibitem[I]{I}  {\em A.~Iarrobino}, Punctual Hilbert Schemes, Mem.\
Amer.\ Math.\ Soc. {\bf 188} (1977).

\bibitem[IY]{IY}  {\em A.~Iarrobino} and {\em J.~Yam\'eogo}, Graded
Ideals in $k[x,y]$ and Partitions: I \& II, preprints, Nice 1992.

\bibitem[K]{K}  {\em A.~D.~King}, Moduli of Representations of Finite
Dimensional Algebras, Quart.\ J.\ Math.\ Oxford, to appear.


\bibitem[L]{L}  {\em J.~Le Potier}, A propos de la construction de
l'espace de modules des faisceaux semi-stables sur le plan projectif,
Bull.\ Soc.\ Math.\ France, to appear.

\bibitem[PS]{PS}  {\em R.~Piene} and {\em M.~Schlessinger}, On the
Hilbert Scheme Compactification of the Space of Twisted Cubics, Amer.\
J.\ Math.\ {\bf 107} (1985), 761--774.

\bibitem[S]{S} {\em J.-P.~Serre}, Faisceaux alg\'ebriques coh\'erents,
Ann.\ of Math.\ {\bf 61} (1955), 197--278.

\bibitem[Y]{Y}  {\em J.~Yam\'eogo}, Fibr\'e canonique de la
vari\'et\'e $G_T$ param\'etrant les id\'eaux homog\`enes de ${\Bbb
C}[[x,y]]$ ayant pour fonction de Hilbert $T$. Preprint, Nice 1993.
\end{thebibliography}
\end{document}